\documentclass{amsart}

\usepackage{amsmath,amsthm,amssymb,graphicx,color,natbib}

\def\@journal{\ }
\newcommand{\eqdef}{\stackrel{d}{=}}

\newcommand{\probconv}{\stackrel{p}{\longrightarrow}}

\newtheorem{thm}{Theorem}[section]
\newtheorem{remark}[thm]{Remark}
\newtheorem{propn}[thm]{Proposition}

\newtheorem{cor}[thm]{Corollary}

\numberwithin{equation}{section}

\begin{document}

\title{Nonsingular Group Actions and Stationary S$\alpha$S Random Fields}


\author{Parthanil Roy}

\address{Department of Statistics and Probability, Michigan State University, East Lansing, MI 48824-1027, USA}

\email{roy@stt.msu.edu}

\thanks{Supported in part by NSF
grant DMS-0303493 and NSF training grant ``Graduate and Postdoctoral
Training in Probability and Its Applications'' at Cornell
University, the RiskLab of the Department of Mathematics, ETH Zurich and a start up grant from the Michigan State University.}



\begin{abstract}
This paper deals with measurable stationary symmetric stable random fields indexed by $\mathbb{R}^d$ and their relationship with the ergodic theory of nonsingular $\mathbb{R}^d$-actions. Based on the phenomenal work of \cite{rosinski:2000}, we establish extensions of some structure results of stationary $S\alpha S$ processes to $S\alpha S$ fields. Depending on the ergodic theoretical nature of the underlying action, we observe different behaviors of the extremes of the field.
\end{abstract}

\maketitle

\section{Introduction}

$\mathbf{X}:=\{X_t\}_{t\in \mathbb{R}^d}$ is called a symmetric
$\alpha$-stable ($S\alpha S$) random field if for all $c_1, c_2,
\ldots, c_k \in \mathbb{R}$ and $t_1,t_2,\ldots,t_k \in
\mathbb{R}^d$, $\sum_{j=1}^k c_j X_{t_j}$ follows a symmetric
$\alpha$-stable distribution. See \cite{samorodnitsky:taqqu:1994} for more information on $S \alpha S$ distributions and processes. In this paper we will further assume that $\{X_t\}_{t\in \mathbb{R}^d}$ is measurable and stationary with $\alpha \in (0,2)$.

The Hopf decomposition of nonsingular flows (see \cite{aaronson:1997}) gives rise to a useful decomposition of stationary $S\alpha S$ processes into two independent components; see \cite{rosinski:1995}. For a general $d>1$, \cite{rosinski:2000} established a similar decomposition of $S\alpha S$ random fields. We show the connection between this work and the conservative-dissipative decomposition of nonsingular $\mathbb{R}^d$-actions. This connection with ergodic theory enables us to study the rate of growth of the partial maxima $\{M_\tau\}_{\tau>0}$ of the random field $X_t$ as $t$ runs over a $d$-dimensional hypercube with an edge length $\tau$ increasing to infinity. This is a straightforward extension of the one-dimensional version of this result available in \cite{samorodnitsky:2004b}. See \cite{samorodnitsky:2004a} and \cite{roy:samorodnitsky:2008} for the discrete parameter case.

This paper is organized as follows. In Section \ref{sec_ergodic_theory_cont}, we develop the theory of nonsingular $\mathbb{R}^d$-actions based on \cite{aaronson:1997} and  \cite{kolodynski:rosinski:2003}. We extend some of the structure results of stationary $S\alpha S$ processes available in \cite{rosinski:1995} to the $d>1$ case in Section \ref{sec_stable_random_fields_cont} and use these results in Section \ref{sec_M_tau_cont} to compute the rate of growth of the partial maxima $M_\tau$ of the field as $\tau$ increases to infinity.

\section{Nonsingular $\mathbf{\mathbb{R}^d}$-actions}\label{sec_ergodic_theory_cont}

In this section we present the theory of nonsingular $\mathbb{R}^d$-actions in parallel to the corresponding discrete-parameter results discussed in Section 2 in \cite{roy:samorodnitsky:2008}. Most of the notions discussed in this section can be found in \cite{aaronson:1997} and \cite{krengel:1985}.

Let $\{\phi_t\}_{t \in \mathbb{R}^d}$ be a nonsingular $\mathbb{R}^d$-action on a $\sigma$-finite standard
measure space $(S, \mathcal{S}, \mu)$. This means that $\{\phi_t\}_{t \in \mathbb{R}^d}$ is a collection of measurable transformations $\phi_t:S\to S$ such that\smallskip\\
\indent (i) $\phi_0(s) = s \mbox{ for all } s\in S$,\smallskip\\
\indent (ii) $\phi_{v+u}(s) = \phi_u\circ\phi_v(s)  \mbox{ for all }  s\in S, u,v\in \mathbb{R}^d$,\smallskip\\
\indent (iii) $(s,u)\mapsto \phi_u(s)$ is measurable map,\smallskip\\
\indent (iv) $\mu \sim \mu \circ \phi_t^{-1}$ for all $t \in \mathbb{R}^d$.\smallskip\\ 
Define lattices $\Gamma_n:=\frac{1}{2^{n}}\mathbb{Z}^d \subseteq \mathbb{R}^d$ for all $n \geq 0$. The following result is a partial extension of Corollary $1.6.5$ in \cite{aaronson:1997} to nonsingular $\mathbb{R}^d$-actions.

\begin{propn} \label{propn_R^d_cont}
Conservative (resp. dissipative) parts of the actions $\{\phi_t\}_{t
\in \Gamma_n}$, $n \geq 0$, are all equal modulo $\mu$.
\end{propn}

\begin{proof} Let $\mathcal{C}_n$ be the conservative part of $\{\phi_t\}_{t \in \Gamma_n}$ for all $n \geq 0$ and $\lambda$ be the Lebesgue measure on $\mathbb{R}^d$. By Theorem A.1 in \cite{kolodynski:rosinski:2003}, there exists a strictly positive measurable function $(t,s)\mapsto w_t(s)$ defined on
$\mathbb{R}^d \times S$, such that
for all $t\in {\mathbb R}^d$,
\[
w_t(s) = \frac{d\mu\circ \phi_t}{d\mu}(s)
\]
for $\mu$-almost all $s \in S$, and for all $t,h \in \mathbb{R}^d$ and for all $s \in S$
\begin{equation}\label{e:w(t+h)}
 w_{t+h}(s) = w_h(s) w_t(\phi_h(s)).
\end{equation}

Let, for all $n \geq 0$, $F_n:=[\mathbf{0},\frac{1}{2^n}\mathbf{1})$, where
$\mathbf{0}=(0,0,\ldots,0),\,\mathbf{1}=(1,1,\ldots,1) \in \mathbb{R}^d$ and for all
$u=(u^{(1)},u^{(2)},\ldots,u^{(d)}),\,v=(v^{(1)},v^{(2)},\ldots,v^{(d)}) \in \mathbb{R}^d$, $[u,v):=\{x \in \mathbb{R}^d:\,u^{(i)}
\leq x^{(i)} < v^{(i)} \mbox{ for all }i=1,2,\ldots,d\}$. Taking $h \in L^1(S,\mu)$, $h>0$, and using \eqref{e:w(t+h)}, we get, for all $s \in
S$ and for all $n\geq 0$,
\begin{align*}
\int_{\mathbb{R}^d} h\circ \phi_t(s)w_t(s)\lambda(dt) 
    &=\sum_{\gamma \in \Gamma_n}\int_{F_n} h\circ \phi_{\gamma+t}(s)w_{\gamma+t}(s)\lambda(dt)\\
    &=\sum_{\gamma \in \Gamma_n} h_n \circ \phi_{\gamma}(s)w_{\gamma}(s),
\end{align*}
where $h_n(s):=\int_{F_n} h\circ \phi_t(s)w_t(s)\lambda(dt) \in L^1(S,\mu)$
by Fubini's theorem. Hence, by Corollary 2.4 in \cite{roy:samorodnitsky:2008}, we get that for all $n \geq 0$,
\[
\mathcal{C}_n =\left\{s \in S:\,\int_{\mathbb{R}^d} h\circ \phi_t(s)w_t(s)\lambda(dt)=\infty\right\} \mbox{ modulo }\mu,
\]
which completes the proof.
\end{proof}

Motivated by Proposition \ref{propn_R^d_cont}, we define the conservative (resp. dissipative) part of $\{\phi_t\}_{t\in \mathbb{R}^d}$ to be $\mathcal{C}_0$
(resp. $\mathcal{D}_0:=S \setminus \mathcal{C}_0$). Then from the
proof of Proposition \ref{propn_R^d_cont}, we get the following continuous parameter analogue of Corollary 2.4 in \cite{roy:samorodnitsky:2008}.
\begin{cor} \label{cor_R^d_cont}
For any $h \in L^1(S,\mu)$, $h>0$, the conservative part of
$\{\phi_t\}_{t\in \mathbb{R}^d}$ is given by
\[
\mathcal{C}=\left\{s \in S:\,\int_{\mathbb{R}^d} h\circ
\phi_t(s)w_t(s)\lambda(dt)=\infty\right\} \mbox{ modulo } \mu,
\]
where $w_t(s)$ is as above.
\end{cor}

\begin{remark} \textnormal{Note that Theorem A.1 in \cite{kolodynski:rosinski:2003} takes care of the measurability issues regarding the Radon Nikodym derivatives very nicely.}
\end{remark}

As in the discrete case, the action $\{\phi_t\}$ is called conservative if $S=\mathcal{C}$ and dissipative if $S=\mathcal{D}$. Recall that nonsingular group actions
$\{\phi_t\}_{t\in \mathbb{R}^d}$ and $\{\psi_t\}_{t\in \mathbb{R}^d}$, defined on standard
measure spaces $(S,\mathcal{S},\mu)$ and $(T,\mathcal{T},\nu)$
resp., are equivalent if there is a Borel isomorphism $\Phi:S \to T$ such that $\nu\sim \mu \circ \Phi^{-1}$ and for
each $t \in \mathbb{R}^d$, $$\psi_t \circ \Phi=\Phi \circ \phi_t $$  $\mu$-almost surely. In light of Corollary \ref{cor_R^d_cont}, we can rephrase Theorem $2.2$ in \cite{rosinski:2000} to obtain Krengel's structure theorem (see \cite{krengel:1969}) for
dissipative nonsingular $\mathbb{R}^d$-actions.
\begin{cor}[\cite{rosinski:2000}]\label{cor_krengel_thm_R^d}
Let $\{\phi_t\}$ be a nonsingular $\mathbb{R}^d$-action on a
$\sigma$-finite standard measure space $(S,\mathcal{S},\mu)$. Then
$\{\phi_t\}$ is dissipative if and only if it is equivalent to the
$\mathbb{R}^d$-action
$\psi_t(w,s):=(w,t+s)
$
defined on $(W \times {\mathbb{R}}^d, \tau \otimes \lambda)$, where $(W,
\mathcal{W},\tau)$ is some $\sigma$-finite standard measure space
and $\lambda$ is the Lebesgue measure on $\mathbb{R}^d$.
\end{cor}

\section{Structure of Stationary S$\alpha$S Random Fields}\label{sec_stable_random_fields_cont}

Suppose $\mathbf{X}=\{X_t\}_{t \in \mathbb{R}^d}$ is a stationary measurable
$S\alpha S$ random field, $0 < \alpha < 2$. Every
measurable minimal representation (this exists by Theorem $2.2$ in \cite{rosinski:1995}) of $\mathbf{X}$ is of the from
\begin{eqnarray}
X_t&\eqdef& \int_{S} f_t(s)M(ds),\;\; t \in \mathbb{R}^d, \nonumber
\end{eqnarray}
where
\begin{eqnarray}
f_t(s)=c_t(s){\left(\frac{d \mu \circ \phi_t}{d
\mu}(s)\right)}^{1/\alpha}f \circ \phi_t(s)\label{integral_repn_stationary_cont}
\end{eqnarray}
for all $t \in \mathbb{R}^d$ and $s \in S$, $M$ is an $S\alpha S$ random measure on some standard Borel
space $(S,\mathcal{S})$ with $\sigma$-finite control measure $\mu$, $f \in L^{\alpha}(S,\mu)$, $\{\phi_t\}_{t \in \mathbb{R}^d}$
is a nonsingular $\mathbb{R}^d$-action on $(S, \mu)$ and $\{c_t\}_{t
\in \mathbb{R}^d}$ is a measurable cocycle for $\{\phi_t\}$ taking
values in $\{-1, +1\}$, i.e., $(t,s) \mapsto c_t(s)$ is a jointly measurable map $\mathbb{R}^d \times S \rightarrow \{-1, +1\}$ such that for all $u,v \in \mathbb{R}^d$, $c_{u+v}(s)=c_v(s) c_u\big(\phi_v(s)\big)$ for $\mu$-a.a. $s
\in S$; see \cite{rosinski:1995} for the $d=1$ case and
\cite{rosinski:2000} for a general $d$. 

Conversely, $\{X_t\}$ defined as above is a stationary measurable $S\alpha S$ random
field. Without loss of generality we can assume that the family
$\{f_t\}$ in \eqref{integral_repn_stationary_cont} satisfies the full support assumption
\begin{eqnarray}
\mbox{Support}\left\{f_t:\,t \in \mathbb{R}^d\right\}=S
\label{condn_full_support_cont}
\end{eqnarray}
and take the Radon-Nikodym derivative in \eqref{integral_repn_stationary_cont} to be equal to $w_t(s)$ defined in Section \ref{sec_ergodic_theory_cont} by virtue of Theorem A.1 in \cite{kolodynski:rosinski:2003}. We first establish that any measurable stationary random
field indexed by $\mathbb{R}^d$ is continuous in probability. The corresponding one-dimensional result was established by \cite{surgailis:rosinski:mandrekar:cambanis:1998} using a result of \cite{cohn:1972}.

\begin{propn} \label{propn_cont_in_prob_cont}
Suppose $\mathbf{X}=\{X_t\}_{t \in \mathbb{R}^d}$ be a measurable
stationary random field. Then $\mathbf{X}$ is continuous in
probability, i.e., for every $t_0 \in \mathbb{R}^d$, $X_t \probconv
X_{t_0}$ whenever $t \rightarrow t_0$.
\end{propn}

\begin{proof} Using a truncation argument we can assume without loss of generality that $\|X_0\|_{2} <
\infty$ where $\|\cdot\|_{2}$ denotes the $L^2$-norm. Define $\{\phi_t\}_{t \in \mathbb{R}^d}$ to be the shift action on the path-space $\Omega$ given by $\phi_t(\omega)(s)=\omega(s+t)$ for all $\omega \in \Omega$. By measurability and stationarity of $\mathbf{X}$, $\{\phi_t\}$ is an $\mathbb{R}^d$-action which preserves the induced probability measure. Using Banach's theorem for Polish groups (see \cite{banach:1932} p.~20) it follows that $t \mapsto X_t$ is $L^2$-continuous (see Section $1.6$ in \cite{aaronson:1997}), which implies the result.
\end{proof}

As in the discrete parameter case, we say that a measurable
stationary $S\alpha S$ random field $\{X_t\}_{t\in \mathbb{R}^d}$ is generated by a nonsingular $\mathbb{R}^d$-action $\{\phi_t\}$ on
$(S, \mu)$ if it has an integral representation of the form
$(\ref{integral_repn_stationary_cont})$ satisfying
$(\ref{condn_full_support_cont})$. The following result, which is the continuous parameter analogue of Proposition $3.1$ in \cite{roy:samorodnitsky:2008}, yields that the classes of measurable stationary $S\alpha S$ random
fields generated by conservative and dissipative actions are
disjoint. The corresponding one-dimensional result is available in Theorem $4.1$ of \cite{rosinski:1995}.

\begin{propn} \label{propn_cons_diss_cont}
Suppose $\{X_t\}_{t\in \mathbb{R}^d}$ is a measurable stationary
$S\alpha S$ random field generated by a nonsingular
$\mathbb{R}^d$-action $\{\phi_t\}$ on $(S, \mu)$ and $\{f_t\}$ is
given by $(\ref{integral_repn_stationary_cont})$. Let
$\mathcal{C}$ and $\mathcal{D}$ be the conservative and dissipative
parts of $\{\phi_t\}$. Then we have 
$$\mathcal{C}=\{s \in S:\int_{\mathbb{R}^d}
|f_t(s)|^{\alpha}\lambda(dt)=\infty
\}$$ 
and 
$$\mathcal{D}=\{s \in S:\int_{\mathbb{R}^d}
|f_t(s)|^{\alpha}\lambda(dt)<\infty \}$$ 
modulo $\mu$. In particular, if a stationary $S\alpha S$ random field $\{X_t\}_{t \in
\mathbb{R}^d}$ is generated by a conservative (dissipative, resp.)
$\mathbb{R}^d$-action, then in any other integral representation of
$\{X_t\}$ of the form $(\ref{integral_repn_stationary_cont})$
satisfying $(\ref{condn_full_support_cont})$, the
$\mathbb{R}^d$-action must be conservative (dissipative, resp.).
\end{propn}

\begin{proof} Let $$h(s):= \sum_{\gamma \in \mathbb{Z}^d} a_\gamma
\int_{\gamma+F_0}|f_t(s)|^\alpha \lambda(dt),$$ where $s \in
S, a_\gamma >0$ for all $\gamma \in \mathbb{Z}^d$ and
$\sum_{\gamma \in \mathbb{Z}^d}a_\gamma=1$. Clearly $h \in L^1(S,
\mu)$ and $h>0$ almost surely. By \eqref{e:w(t+h)} and the translation invariance of $\lambda$,
$$
\sum_{\beta \in \mathbb{Z}^d} h\circ\phi_\beta(s)w_\beta(s)
=\int_{\mathbb{R}^d}|f_t(s)|^{\alpha} \lambda(dt)
$$
for all $s \in S$. Hence, by Corollary 2.4 in \cite{roy:samorodnitsky:2008}, we get
\begin{align*}
\mathcal{C}=\mathcal{C}_0&=\Big\{s \in S:\,\sum_{\beta \in \mathbb{Z}^d} h\circ\phi_\beta(s)w_\beta(s)=\infty\Big\}\\
&=\Big\{s \in S:\,\int_{\mathbb{R}^d} |f_t(s)|^{\alpha}\lambda(dt)=\infty\Big\}\;\;\; \mbox{ modulo }\mu.
\end{align*}
This completes the proof of the first part.

The second part follows by an argument parallel to the one in the proof of Theorem 4.1 in \cite{rosinski:1995}.
\end{proof}

The following corollary is the continuous parameter analogue of Corollary $3.2$ of \cite{roy:samorodnitsky:2008}. The corresponding one-dimensional result is available in Corollary $4.2$ of \cite{rosinski:1995} and the same proof works in the $d$-dimensional case.

\begin{cor} \label{cor_cons_diss_cont} The measurable stationary $S\alpha S$
random field $\{X_t\}_{t \in \mathbb{R}^d}$ is generated by a
conservative (dissipative, resp.) $\mathbb{R}^d$-action if and only
if for any (equivalently, some) measurable representation
$\{f_t\}_{t \in \mathbb{R}^d}$ of $\{X_t\}$ satisfying
$(\ref{condn_full_support_cont})$, the integral $\int_{\mathbb{R}^d} |f_t(s)|^{\alpha}d\lambda(t)$
is infinite (finite, resp) $\mu$-almost surely.
\end{cor}

Recall that \cite{surgailis:rosinski:mandrekar:cambanis:1993} defined $\mathbf{X}$ to be a stable mixed
moving average if
\begin{eqnarray}
\mathbf{X} \eqdef \left\{\int_{W \times
{\mathbb{R}}^d}f(v,t+s)\,M(dv,ds)\right\}_{t \in {\mathbb{R}}^d}\,,
\label{defn_mixed_moving_avg_cont}
\end{eqnarray}
where $f \in L^{\alpha}(W \times {\mathbb{R}}^d, \nu \otimes
\lambda)$, $\lambda$ is the Lebesgue measure on ${\mathbb{R}}^d$,
$\nu$ is a $\sigma$-finite measure on a standard Borel space $(W,
\mathcal{W})$, and the control measure $\mu$ of $M$ equals $\nu
\otimes \lambda$. The following result gives three equivalent
characterizations of stationary $S\alpha S$ random fields generated
by dissipative actions.

\begin{thm} \label{thm_mixed_moving_avg_cont} Suppose $\{X_t\}_{t \in
{\mathbb{R}}^d}$ is a measurable stationary $S\alpha S$ random field. Then the
following are equivalent:
\begin{enumerate}
\item $\{X_t\}$ is generated by a dissipative $\mathbb{R}^d$-action.

\item For any measurable representation $\{f_t\}$ of $\{X_t\}$ we have,
\[
\int_{\mathbb{R}^d} |f_t(s)|^{\alpha}<\infty \mbox{ for }
\mu\mbox{-a.a. }s.
\]
\item $\{X_t\}$ is a mixed moving average.

\item $\{X_t\}_{t \in \Gamma_n}$ is a mixed moving average for some
(all) $n \geq 1$.
\end{enumerate}

\end{thm}

\begin{proof} $(1)$ and $(2)$ are equivalent by Corollary
$\ref{cor_cons_diss_cont}$, $(2)$ and $(3)$ are equivalent by
Theorem $2.1$ of \cite{rosinski:2000}. $(1)$ and $(4)$ are equivalent by
Theorem $3.3$ in \cite{roy:samorodnitsky:2008} and Proposition \ref{propn_R^d_cont}.
\end{proof}

Therefore, in order to verify that $\mathbf{X}$ is a mixed moving average, it is enough to
verify it on a discrete skeleton (e.g., $\{X_t\}_{t \in \mathbb{Z}^d}$) of the random field. Theorem
$\ref{thm_mixed_moving_avg_cont}$ allows us to describe the decomposition of a stationary $S \alpha S$ random field given in
Theorem $3.7$ of \cite{rosinski:2000} in terms of the
ergodic-theoretical properties of  nonsingular
$\mathbb{R}^d$-actions generating the field. See Corollary $3.4$ in \cite{roy:samorodnitsky:2008} for the corresponding discrete parameter result.

\begin{cor} \label{cor_decomp_cont}
A stationary $S \alpha S$ random field $\mathbf{X}$ has a unique in
law decomposition
\begin{equation}
X_t \eqdef X^{\mathcal{C}}_t+X^{\mathcal{D}}_t \label{cons_diss_decomp_cont}
\end{equation}
where $\mathbf{X^{\mathcal{C}}}$ and $\mathbf{X^{\mathcal{D}}}$ are
two independent stationary $S\alpha S$ random fields such that
$\mathbf{X^{\mathcal{D}}}$ is a mixed moving average, and
$\mathbf{X^{\mathcal{C}}}$ is generated by a conservative action.
\end{cor}

\section{A Note on the Extreme Values} \label{sec_M_tau_cont}

The extreme values of $\{X_t\}$ are expected to grow at a slower rate if $\{X_t\}$ is generated by a
conservative action because of longer memory; see, for example, \cite{samorodnitsky:2004a}, \cite{samorodnitsky:2004b} and \cite{roy:samorodnitsky:2008}. This can be formally proved provided $\mathbf{X}=\{X_t\}_{t\in\mathbb{R}^d}$ is assumed to be locally bounded apart from being stationary and measurable. If further $\mathbf{X}$ is
separable then
\begin{equation}
M_\tau=\sup_{\mathbf{0} \leq s \leq \tau\mathbf{1}} |X_s|,\;\;\;\tau
> 0, \label{defn_of_M_tau_crude_cont}
\end{equation}
is a well-defined finite-valued stochastic process. Here $u =
(u^{(1)}, \ldots , u^{(d)})$ $\leq v =(v^{(1)},
\ldots, v^{(d)})$ means $u^{(i)} \leq v^{(i)}$ for all
$i=1,2,\ldots,d$ and $\mathbf{1}:=(1,1,\ldots,1)$,
$\mathbf{0}:=(0,0,\ldots,0)$.

Since $\mathbf{X}$ is stationary and
measurable, it is continuous in probability by Proposition
\ref{propn_cont_in_prob_cont}. Therefore, as in the one-dimensional case in \cite{samorodnitsky:2004b}, taking its separable
version the above maxima process can be defined by
\begin{equation}
M_\tau= \sup_{s \in [\mathbf{0},\,\tau \mathbf{1}]\cap \Gamma}
|X_s|,\;\;\;\tau > 0, \nonumber
\end{equation}
where $\Gamma:=\bigcup_{n=1}^\infty \Gamma_n=\bigcup_{n=1}^\infty
\frac{1}{2^n}\mathbb{Z}^d$ and $[u,v]:=\{s \in \mathbb{R}^d:\,u \leq
s \leq v\}$. This will avoid the usual measurability problems of
the uncountable maximum \eqref{defn_of_M_tau_crude_cont}. The next result is the
continuous parameter extension of Theorem
$4.3$ in \cite{roy:samorodnitsky:2008}. It follows by the exact same argument as in the one-dimensional
version of this result (Theorem $2.2$ in
\cite{samorodnitsky:2004b}) based on Theorem \ref{thm_mixed_moving_avg_cont} and Corollary \ref{cor_decomp_cont}.

\begin{thm} \label{thm_M_tau_rate_cont}
Let $\mathbf{X}=\{X_t\}_{t \in \mathbb{R}^d}$ be a stationary, locally bounded
$S\alpha S$ random field, where $0 < \alpha < 2$. \vspace{0.07in}\\
\noindent (i) Suppose that $\mathbf{X}$ is not generated by a
conservative action (i.e. the component $\mathbf{X}^{\mathcal{D}}$ in
$(\ref{cons_diss_decomp_cont})$ generated by the dissipative part is
nonzero). Then
\begin{eqnarray}
\frac{1}{\tau^{d/\alpha}}M_\tau  \Rightarrow  C^{\,1/\alpha}_\alpha K_X
Z_\alpha \nonumber
\end{eqnarray}
as $\tau \rightarrow \infty$, where $$K_X ={\left(\int_W {(g(v))}^\alpha \nu(dv) \right)}^{1/\alpha},$$ with $$g(v):=\sup_{s \in \Gamma}|f(v,s)|,\;\;v \in W,$$ for any representation of
$\mathbf{X}^{\mathcal{D}}$ in the mixed moving average form
$(\ref{defn_mixed_moving_avg_cont})$, $C_\alpha$ is the stable tail
constant (see (1.2.9) in \cite{samorodnitsky:taqqu:1994})
and $Z_\alpha$ is the standard Fr\'{e}chet-type extreme value random
variable with distribution $$P(Z_\alpha \leq z)=e^{-z^{-\alpha}}$$ for $z>0$.\\

\noindent (ii) Suppose that $\mathbf{X}$ is generated by a
conservative $\mathbb{R}^d$-action. Then
\begin{eqnarray}
\frac{1}{\tau^{d/\alpha}}M_\tau \probconv 0 \nonumber
\end{eqnarray}
as $\tau \rightarrow \infty$. Furthermore, defining $$b_\tau:=\left(\int_S \sup_{t \in [\mathbf{0},\tau \mathbf{1}]\cap \Gamma}
|f_t(s)|^\alpha \mu(ds) \right)^{1/\alpha},$$ we have that $\left\{{c_\tau}^{-1}M_\tau:\,\tau>0 \right\}$ is not tight for any
positive $c_\tau=o(b_\tau)$. If, for some $\theta >0$ and $c >0$,
\begin{eqnarray}
b_\tau  \geq  c\tau^\theta \;\;\;\;\text{ for all $\tau$ large enough,} \label{inequality_cont}
\end{eqnarray}
then $\left\{{b_\tau}^{-1}M_\tau:\,\tau>0\right\}$ is tight. Finally, for $\tau >0$, let $\eta_\tau$ be a
probability measure on $(S, \mathcal{S})$ with $$\frac{d \eta_\tau}{d \mu}(s)=b_\tau^{-\alpha} \sup_{t \in [\mathbf{0},
\tau\mathbf{1}]\cap \Gamma} |f_t(s)|^\alpha$$ for all $s \in S$ and let $U_j^{(\tau)},\, j=1,2$ be independent $S$-valued random
variables with common law $\eta_\tau$. Suppose that
$(\ref{inequality_cont})$ holds and for any $\epsilon
> 0$,
\begin{eqnarray}
&&P\bigg(\mbox{for some }t \in [\mathbf{0},
\tau\mathbf{1}]\cap \Gamma,\nonumber\\
&&\;\;\;\;\;\;\;\;\;\frac{|f_t(U_j^{(\tau)})|}{\sup_{u \in [\mathbf{0},
\tau\mathbf{1}]\cap \Gamma} |f_u(U_j^{(\tau)})|} > \epsilon,\; j=1,2 \bigg)\rightarrow
0 \label{suff_condn_cont}
\end{eqnarray}
as $\tau \rightarrow \infty$. Then
\begin{eqnarray}
\frac{1}{b_\tau}M_\tau  \Rightarrow  C^{\,1/\alpha}_\alpha Z_\alpha
\nonumber
\end{eqnarray}
as $\tau \rightarrow \infty$. A sufficient condition for
$(\ref{suff_condn_cont})$ is $\lim_{\tau \rightarrow \infty} \tau^{-d/2 \alpha}b_\tau = \infty$.
\end{thm}

Theorem \ref{thm_M_tau_rate_cont} gives the exact rate of growth of the maxima only when the underlying group action is not conservative. In the conservative case, the exact rate depends on the group action as well as on the kernel (see the examples in \cite{samorodnitsky:2004a}, \cite{samorodnitsky:2004b} and \cite{roy:samorodnitsky:2008}). For instance, by an obvious extension of Example $6.1$ in \cite{roy:samorodnitsky:2008} to the continous parameter case, it can be observed that the maxima can grow both polynomially as well as logarithmically and it can even converge to a nonextreme value limit after proper normalization.

In the discrete parameter case, depending on the group theoritic properties of the underlying action, a better estimate of this rate is given in \cite{roy:samorodnitsky:2008}; see also \cite{roy:2007}. This connection with abelian group theory is still an open problem in the continuous parameter case and hence needs to be investigated. Two more open problems related to this work are extensions of the results of \cite{samorodnitsky:2005a} and \cite{roy:2007a} to the $d$-dimensional case.\\

\noindent \textbf{Acknowledgment. }The author is thankful to Gennady Samorodnitsky for many useful discussions, to Paul Embrechts for the support during his stay at RiskLab and to the anonymous referee for his/her valuable suggestions all of which contributed significantly to this work.

\end{document}